\documentclass[11pt]{article}
\usepackage{graphics,amsmath,amssymb}
\usepackage{latexsym}
\usepackage{epsfig}
\usepackage{hyperref}
\usepackage{xcolor}
\usepackage{tikz}
\usetikzlibrary{calc,patterns,decorations.pathmorphing,decorations.markings,arrows}
\usepackage{array}   

\def\a{\alpha}
\def\be{\beta}

\def\epsilon{\varepsilon}
\def\ga{\gamma}

\def\la{\lambda}
\def\La{\Lambda}
\def\phi{\varphi}
\def\si{\sigma}

\def\om{\omega}

\newtheorem{theorem}{Theorem}[section]
\newtheorem{lemma}[theorem]{Lemma}

\newtheorem{corollary}[theorem]{Corollary}

\newtheorem{proposition}[theorem]{Proposition}
\newtheorem{remark}[theorem]{Remark}

\def\Z{{\mathbb Z}}
\def\N{{\mathbb N}}
\def\C{{\mathbb C}}
\def\R{{\mathbb R}}
\def\d{\,\hbox{\rm d}}

\newenvironment{Proof}{\removelastskip\par\medskip
\noindent{\em Proof.} \rm}{\penalty-20\null$\square$\par\medbreak}

\newenvironment{Proofy}{\removelastskip\par\medskip
\noindent{\em Proof of Theorem~\ref{th:inv.ingham}.} \rm}{\penalty-20\null$\square$\par\medbreak}

\newenvironment{Proofv}{\removelastskip\par\medskip
\noindent{\em Proof of Theorem~\ref{th:dir.ingham}.} \rm}{\penalty-20\null$\square$\par\medbreak}

\newenvironment{Proof1}{\removelastskip\par\medskip
\noindent{\em Proof of Theorem~\ref{th:diringham}.} \rm}{\penalty-20\null$\square$\par\medbreak}

\newenvironment{Proof2}{\removelastskip\par\medskip
\noindent{\em Proof of Proposition~\ref{pr:haraux-inv}.} \rm}{\penalty-20\null$\square$\par\medbreak}

\newenvironment{Proof3}{\removelastskip\par\medskip
\noindent{\em Proof of Theorem~\ref{th:inv.ingham1}.} \rm}{\penalty-20\null$\square$\par\medbreak}

\title{\bf Viscoelastic aspects of  \\ glass  relaxation models}

\author{Paola Loreti
\thanks{Dipartimento di Scienze di Base e Applicate per l'Ingegneria,
Sapienza Universit\`a di Roma,
Via Antonio Scarpa 16, 00161 Roma (Italy); e-mail: $<$paola.loreti@uniroma1.it$>$ }
\and Daniela Sforza
\thanks{Dipartimento di Scienze di Base e Applicate per l'Ingegneria, 
Sapienza Universit\`a di Roma,
Via Antonio Scarpa 16, 00161 Roma (Italy); e-mail: $<$daniela.sforza@uniroma1.it$>$ }}

\begin{document}
\date{}

\maketitle

\begin{abstract}
We take advantage of the approximation of the stretched exponential function
with a general Prony series in glass relaxation to give some results about the spectral analysis for the equation of viscoelasticity.
Moreover, in the case of the Burgers model we carry out a complete investigation that leads to the representation of the solution.
\end{abstract}


\
\section{Introduction}
The field of glass science has a  long  and significant history. See as main references \cite{MZ, MPVP}.
The interest devoted to glass materials is mainly concerned with high-tech applications regarding the best possible performances for computer displays, see
\cite{ZM}. For some glass  relaxation models the stretched exponential function, obtained by inserting a fractional power into the exponential, has been proposed as stress relaxation modulus
%
\begin{equation*} 
E(t)= e^{-t^\beta},
\end{equation*} 
where $\beta$ is the stretching exponent (a real number between $0$ and $1$).

The connection between stretched exponentials and glass relaxation 
goes back to 1854 \cite{Ko,CCM}. From these seminal papers a long study was done.
 Here we refer to \cite{MM} for a detailed description. When subject to shaping temperatures, glass shows viscoelasticity in deformation. 
Starting from \cite{MM} we consider the viscoelastic approach developed in the book \cite{RHN}
to  understand the problem 
\begin{equation}\label{eq:problem0I}
u_{tt} (t,x,y)=
 \triangle u(t,x,y)-\beta\int_{0}^t\frac {e^{-(t-s)^\beta}}{ (t-s)^{1-\beta}}\triangle u(s,x,y)ds\,,\quad
t\ge 0\,,\,\, (x,y)\in\Omega,
\end{equation} 
where $\triangle$ denotes the Laplace operator in a  disk 
$\Omega$ of radius $R$ in $\R^2$. 
The motivations for considering a disk for the set $\Omega$ are given  by \cite{ZYMK}.
For other references related to viscoelasticity see \cite{CN, D1, D2, DE, LPC}.

However \eqref{eq:problem0I} is an integro-differential equation with a memory kernel having an integrable singularity in $t=0$. Such
problem is difficult to handle, in fact to our knowledge there are no results in literature about spectral analysis for \eqref {eq:problem0I}. 
Motivated by the the goodness of the approximation of the stretched exponential function
with Prony series, see \cite{MM} 
\begin{equation*}
e^{-t^{\beta}}\approx \sum_{i=1}^N  s_ie^{-r_i t}
\qquad
\Big(s_i, r_i>0\,, \quad\sum_{i=1}^N  s_i=1
\Big)\,,
\end{equation*}
in this paper we  consider  the integro-differential equation
\begin{equation}\label{eq:Mproblem0I}
u_{tt}=
 \triangle u-\sum_{i=1}^N b_i \int_0^t e^{- r_i(t-s)} \triangle u(s)ds\,,\quad
t\ge 0\,,\,\, (x,y)\in\Omega.
\end{equation} 
Here we will show a first result on the spectral analysis for equation \eqref{eq:Mproblem0I}. Indeed, we will prove that
for any  Prony series the equation \eqref{eq:Mproblem0I} has always a null eigenvalue and the sum of all its eigenvalues is given by 
$-\sum_{i=1}^{N}r_i$, being $r_i$  
the exponents of the Prony series. 
As expected result about the spectral analysis we presume that for any $N$ the
principal two branches of complex eigenvalues have  imaginary part going to $\infty$ and bounded real part as $\lambda\to\infty$. Moreover, due to the relaxation,  it is very likely that, with the exception of the null eigenvalue, there are also $N-1$ 
branches of real eigenvalues having a negative accumulation point. 

In order to obtain more precise results, simplification of the equation is necessary. Mechanical models involving springs and dashpots are used to explain the creep and the stress relaxation of viscoelastic deformations.
Among various mechanical models, Burgers model is a typical model which combine a series of elements with springs and dashpots and describe the case in which a Maxwell and a Kelvin-Voigt model are connected in series. To consider the Burgers model is, in fact, a simplification, because the corresponding equation of the viscoelasticity has  as memory kernel a Prony series with $N=2$. For the Burgers model we are able to perform a complete and detailed spectral analysis. 
In particular, we give asymptotic behaviour of all eigenvalues that allows us to represent the solution of the integro-differential equation as a Fourier series.


\section{The Prony series representation of stretched exponential relaxation}\label{se:Mauro}

In a material with memory the stress depends on the entire temporal history of the strain. The linearized constitutive relation for small deformations given in 1874 by Boltzmann \cite{Boltz} leads to the following integro-differential equation
\begin{equation}\label{eq:problem0}
u_{tt} (t,x,y)=
c^2 \triangle u(t,x,y)-\int_{0}^t\ m(t-s) \triangle u(s,x,y)ds\,,\quad
t\ge 0\,,\,\, (x,y)\in\Omega\,,
\end{equation}
where $\triangle$ represents the Laplace operator in a  disk 
$\Omega$ of radius $R$ in $\R^2$. 
Here the constant
\begin{equation}\label{eq:c2}
c^2:=\alpha+\int^{\infty}_0 \ m(s)ds
\qquad 
(\alpha\ge0)
\end{equation}
measures the instantaneous response of stress to strain and is called the instantaneous stress modulus and 
the integral  kernel $m(t)$ can be deduced by means of a so-called stress relaxation test, see \cite{RHN}.
Indeed, if we set the strain $\varepsilon=0$ for $t<0$ and $\varepsilon=\varepsilon_0$ for $t>0$, the stress $\sigma(t)$ for $t>0$ is given by
\begin{equation*}
\sigma(t)=\Big(\alpha+\int^{\infty}_t \ m(s)ds\Big)\varepsilon_0
\,.
\end{equation*}
By measuring the stress, since $\varepsilon_0$ is a constant value
 one obtains the stress relaxation modulus $E(t)$, that is defined as
 \begin{equation}\label{const-rel}
E(t):=\alpha+\int^{\infty}_t \ m(s)ds
\,.
\end{equation} 
From the above formula we derive $E(0)=c^2$ and  the expression of  the memory kernel in terms of the relaxation function $E(t)$, that is
\begin{equation}\label{eq:mE}
m(t)= -E'(t)
\,.
\end{equation}
In the applications for glass models the stress relaxation modulus can be taken as the
stretched exponential function  
\begin{equation}
E(t)=e^{-t^{\beta}}, \quad 0<\beta\le1\,,
\end{equation}
but the above definition leads to introduce singular memory kernels. Indeed, thanks to
\eqref{eq:mE} we have
\begin{equation*}
m(t)= \frac{\beta}{t^{1-\beta}}\ e^{-t^{\beta}}
\,.
\end{equation*}
To overcome the problems deriving from singular  kernels, a mathematical convenient way 
is to represent the stretched exponential function as a Prony series (see \cite{MM} and references therein) i.e. as a discrete sum of simple exponential terms:
\begin{equation}
e^{-t^{\beta}}\approx \sum_{i=1}^N  s_ie^{-r_i t}
\qquad
s_i, r_i>0\,, N\in\N,
\end{equation}
with the weighting factors $s_i$ satisfying
\begin{equation}\label{eq:si}
\sum_{i=1}^N  s_i=1
\,.
\end{equation}
So the relaxation function is 
\begin{equation}
E(t)=\sum_{i=1}^N  s_ie^{-r_i t}
\,,
\end{equation}
whence  the memory kernel is given by
\begin{equation}\label{eq:Mm}
m(t)=-E'(t)=\sum_{i=1}^N s_i r_ie^{-r_i t}=\sum_{i=1}^N b_ie^{-r_i t}
\,,
\end{equation}
that is $b_i=s_i r_i$, $i=1,\dotsc,N$, and the instantaneous stress modulus is
\begin{equation}\label{eq:Mc}
c^2=E(0)=\sum_{i=1}^N  s_i=1
\,.
\end{equation}
Now, taking into account \eqref{eq:Mm} and \eqref{eq:Mc}  we can write  the integro-differential equation \eqref{eq:problem0} in the form
\begin{equation}\label{eq:Mproblem0}
u_{tt}=
 \triangle u-\sum_{i=1}^N b_i \int_0^t e^{- r_i(t-s)} \triangle u(s)ds\,,\quad
t\ge 0\,,
\end{equation} 
with 
\begin{equation}\label{eq:biri}
\sum_{i=1}^N  \frac{b_i}{r_i}=1
\,,
\end{equation}
in virtue of \eqref{eq:si}.
Our goal is to show that for any Prony series the equation \eqref{eq:Mproblem0} has always a null eigenvalue and the sum of all its eigenvalues is given by 
$-\sum_{i=1}^{N}r_i$, being $r_i$  
the exponents of the Prony series. 
First, in \eqref{eq:Mproblem0} we replace the operator $-\triangle$ with its generic eigenvalue $\lambda>0$, that is
\begin{equation}\label{eq:Meig}
u{''}=-\la u+\la \sum_{i=1}^N b_i \int_0^t e^{- r_i(t-s)}u(s) ds\,,
\qquad t\ge 0.
\end{equation}
To write the equation for the eigenvalues, we introduce the variables
\begin{equation*}
v=u',
\quad
w_i=e^{- r_i t}*u,
\quad
i=1,\dotsc,N,
\end{equation*}
and note that the integro-differential equation \eqref{eq:Meig} is equivalent to the following system of first order differential equations
\begin{equation*}
\begin{cases}
u'=v   
\\
\displaystyle v'=-\la u+\la \sum_{i=1}^N b_i w_i
\\
w_1'=u-r_1w_1
\\
w_2'=u-r_2w_2
\\
\dots\dots\dots\dots\dots
\\
w_N'=u-r_Nw_N
\end{cases}
\end{equation*}
The $(N+2)\times(N+2)-$matrix  of the  system is given by
\begin{equation}\label{eq:Nmatrix}
A_N=\begin{pmatrix}
0 & 1 & 0 & 0 & 0 & \dots & 0
\\
-\lambda  & 0 & \lambda b_1 & \lambda b_2 & \lambda b_3  & \dots & \lambda b_N
\\
 1 & 0 &- r_1 & 0 & 0 & \dots & 0
 \\
  1 & 0  & 0 &- r_2 & 0 & \dots & 0
  \\
  \dots & \dots  & \dots & \dots & \dots & \dots & \dots
 \\
  1 & 0  & 0 & \dots & \dots & 0 & - r_N
\end{pmatrix}
\end{equation}
The determinant $\big|A_N-zI\big|$ is a $(N+2)-$polynomial in the variable $z$, precisely
\begin{equation}
\big|A_N-zI\big|=\begin{vmatrix}
-z & 1 & 0 & 0 & 0 & \dots & 0
\\
-\lambda  & -z & \lambda b_1 & \lambda b_2  & \dots & \lambda b_{N-1}  & \lambda b_N
\\
 1 & 0 &- r_1-z & 0 & 0 & \dots & 0
 \\
  1 & 0  & 0 &- r_2-z & 0 & \dots & 0
  \\
  \dots & \dots  & \dots & \dots & \dots & \dots & \dots
  \\
  1 & 0  & 0 & \dots & 0 & - r_{N-1}-z & 0 
\\
  1 & 0  & 0 & \dots & \dots & 0 & - r_N-z
\end{vmatrix}
\end{equation}
By solving the determinant according to the last column, we obtain the following recursive formula
\begin{equation*}
\big|A_N-zI\big|
=(-1)^{N}\lambda b_N
\begin{vmatrix}
-z & 1 & 0 & 0 & \dots & 0 
\\
 1 & 0 &- r_1-z & 0 & 0 & \dots 
 \\
  1 & 0  & 0 &- r_2-z & 0 & \dots 
  \\
  \dots & \dots  & \dots & \dots & \dots & \dots 
  \\
  1 & 0  & 0 & \dots & 0 & - r_{N-1}-z  
\\
  1 & 0  & 0 & \dots & \dots & 0 
\end{vmatrix}
-(r_N+z)\big|A_{N-1}-zI\big|
\,,
\end{equation*}
where $A_{N-1}$ is the $(N+1)-$matrix corresponding to the Prony series $\sum_{i=1}^{N-1} b_ie^{-r_i t}$.
Since
\begin{equation*}
\begin{split}
\begin{vmatrix}
-z & 1 & 0 & 0 & \dots & 0 
\\
 1 & 0 &- r_1-z & 0 & 0 & \dots 
 \\
  1 & 0  & 0 &- r_2-z & 0 & \dots 
  \\
  \dots & \dots  & \dots & \dots & \dots & \dots 
  \\
  1 & 0  & 0 & \dots & 0 & - r_{N-1}-z  
\\
  1 & 0  & 0 & \dots & \dots & 0 
\end{vmatrix}
=&(-1)^{N}
\begin{vmatrix}
 1 & 0 & 0 & \dots & 0 
\\
  0 &- r_1-z & 0 & 0 & \dots 
 \\
   0  & 0 &- r_2-z & 0 & \dots 
  \\
  \dots  & \dots & \dots & \dots & \dots 
  \\
   0  & 0 & \dots & 0 & - r_{N-1}-z  
\end{vmatrix}
\\
\\
=&-(r_1+z)(r_2+z)\dotsb(r_{N-1}+z)
\,,
\end{split}
\end{equation*}
we have
\begin{equation}\label{eq:detN}
\big|A_N-zI\big|
=(-1)^{N+1}\lambda b_N
(r_1+z)(r_2+z)\dotsm(r_{N-1}+z)
-(r_N+z)\big|A_{N-1}-zI\big|
\,.
\end{equation}
We will show by induction that for any $N\ge2$ the inhomogeneous term of the polynomial $\big|A_{N}-zI\big|$ is given by
\begin{equation}\label{eq:inhN}
(-1)^{N+1}\lambda
\big(b_1r_2\dotsm r_{N}+b_2r_1r_3\dotsm r_{N}+ \dots\ +b_{N}r_1\dotsm r_{N-1}-r_1\dotsm r_{N}\big)\,.
\end{equation}
If $N=2$ from \eqref{eq:detN} it follows
\begin{equation*}
\big|A_2-zI\big|
=-\lambda b_2
(r_1+z)
-(r_2+z)\big|A_{1}-zI\big|
\,.
\end{equation*}
Since
\begin{equation}\label{eq:detN1}
\big|A_1-zI\big|
=-z^3-r_1 z^2-\lambda z+\lambda(b_1-r_1)
\,,
\end{equation}
we have
\begin{equation*}
\big|A_2-zI\big|
=z^{4}+(r_1+r_2) z^{3}+(\la+r_1r_2) z^{2}+\la (r_1+r_2- b_1- b_2)z-\la (b_1r_2+b_2r_1-r_1r_2)\,,
\end{equation*}
and hence the formula \eqref{eq:inhN} is satisfied.
For an arbitrary $N$ we assume that the inhomogeneous term of the polynomial $\big|A_{N-1}-zI\big|$ is given by
\begin{equation*}
(-1)^{N}\lambda\big(b_1r_2\dotsm r_{N-1}+b_2r_1r_3\dotsm r_{N-1}+\dots\ +b_{N-1}r_1\dotsm r_{N-2}-r_1\dotsm r_{N-1}\big)
\,.
\end{equation*}
Taking into account of the previous formula we get that the inhomogeneous term of the polynomial $\big|A_{N}-zI\big|$ is given by
\begin{multline*}
(-1)^{N+1}\lambda b_N
r_1r_2\dotsm r_{N-1}
-r_N (-1)^{N}\lambda\big(b_1r_2\dotsm r_{N-1}+b_2r_1r_3\dotsm r_{N-1}+ \dots\ +b_{N-1}r_1\dotsm r_{N-2}-r_1\dotsm r_{N-1}\big)
\\
=(-1)^{N+1}\lambda
\big(b_1r_2\dotsm r_{N}+b_2r_1r_3\dotsm r_{N}+ \dots\ +b_{N}r_1\dotsm r_{N-1}-r_1\dotsm r_{N}\big)
\end{multline*}
that is formula \eqref{eq:inhN}, and hence our statement  holds true for any $N$.

From \eqref{eq:inhN} it follows that the equation \eqref{eq:Mproblem0} has always a null eigenvalue. Indeed, in virtue of \eqref{eq:biri} we have
\begin{equation*}
b_1r_2\dotsm r_{N}+b_2r_1r_3\dotsm r_{N}+ \dots\ +b_{N}r_1\dotsm r_{N-1}=r_1\dotsm r_{N}
\,,
\end{equation*}
whence, taking into account \eqref{eq:inhN}, we have the equation of the eigenvalues 
\begin{equation*}
\big|A_N-zI\big|=0
\end{equation*}
has null inhomogeneous term. So, the previous equation  has always the solution $z=0$.

Now, again by induction we will show that for $N\ge1$ the term $z^{N+1}$ of the polynomial $\big|A_{N}-zI\big|$ is 
\begin{equation}\label{eq:zN+1}
(-1)^{N}\big(r_1+r_2+\dotsb + r_{N}\big)z^{N+1}
\,.
\end{equation}
For $N=1$ our assertion follows from
\eqref{eq:detN1}.
In addition, if we assume that the term $z^{N}$ of the polynomial $\big|A_{N-1}-zI\big|$ is 
\begin{equation*}
(-1)^{N-1}\big(r_1+r_2+\dotsm + r_{N-1}\big)z^{N}
\,,
\end{equation*}
thanks to \eqref{eq:detN} we get that  term $z^{N+1}$ of the polynomial $\big|A_{N}-zI\big|$ is given by
\begin{equation*}
-r_N(-1)^{N-1}z^{N+1}-z(-1)^{N-1}\big(r_1+r_2+\dotsb + r_{N-1}\big)z^{N}
=(-1)^{N}\big(r_1+r_2+\dotsb + r_{N}\big)z^{N+1},
\end{equation*}
that is \eqref{eq:zN+1}. 

Finally, recalling that the sum of the zeros of a $(N+2)$-degree polynomial is given by minus the $(N+1)$-degree coefficient, from \eqref{eq:zN+1}  we deduce that the sum of the eigenvalues of the equation \eqref{eq:Mproblem0} is given by $-\sum_{i=1}^{N}r_i$, that is, it depends only on  
the exponents $r_i$ of the Prony series. 

One can perform numerical simulations by means of \eqref{eq:Nmatrix}. Indeed, for some sets of values of $b_i$ and $r_i$ satisfying the condition $\sum_{i=1}^N  \frac{b_i}{r_i}=1$ it is possible to obtain the corresponding expression of the eigenvalues  $z_i$ as in Tables \ref{table:1} and \ref{table:2}.
Such numerical simulations show that the principal two branches of complex eigenvalues have  imaginary part going to $\infty$ and bounded real part as $\lambda\to\infty$. Due to the relaxation, with the exception of the null eigenvalue, there are also $N-1$ branches of real eigenvalues having a negative accumulation point.  It remains an open problem to show such behaviour for any $N$ from a theoretical point of view.

\begin{table}[h!]
\begin{center}
\begin{tabular}{ | m{0.5cm} | m{7.5em}| m{0.01cm} | m{10em}| m{0.01cm} | m{12em}|  } 
\hline
$N$ &$3$\hfill & &$4$ & &$5$  \\
\hline
$b_i$& 1\hfill  2\hfill  3&  & 1\hfill  2\hfill  3\phantom{1}\hfill  4\phantom{1} 
& & 1\hfill  2\phantom{1}\hfill  3\phantom{1}\hfill  4\phantom{1}\hfill  5\phantom{1}\\ 
\hline
$r_i$ & 3\hfill  6\hfill  9&  & 4\hfill  8\hfill  12\hfill 16 
& & 5\hfill  10\hfill  15\hfill  20\hfill 25\\ 
\hline
$z_i$ & $0$, $-4.40$, $-7.95$, 
$ -2.82 + 9.28i$,  
$-2.82 - 9.28i$  & 
& $0$, $-5.94$, $-10.68$, $-15.13$,
$ -4.11 + 8.13i$,  
$-4.11 - 8.13i$ &
& $0$, $-7.84$, $-13.77$, $-19.17$, $-24.40$,
$-4.90 + 6.61i$,
$-4.90 - 6.61i$
\\ 
\hline
\end{tabular}
\end{center}
\caption{Eigenvalues for $\lambda=100$}
\label{table:1}
\end{table}

\begin{table}[h!]
\begin{center}
\begin{tabular}{ | m{0.5cm} | m{7.5em}| m{0.01cm} | m{10em}| m{0.01cm} | m{12em}|  } 
\hline
$N$ &$3$\hfill & &$4$ & &$5$  \\
\hline
$b_i$& 1\hfill  2\hfill  3&  & 1\hfill  2\hfill  3\phantom{1}\hfill  4\phantom{1} 
& & 1\hfill  2\phantom{1}\hfill  3\phantom{1}\hfill  4\phantom{1}\hfill  5\phantom{1}\\ 
\hline
$r_i$ & 3\hfill  6\hfill  9&  & 4\hfill  8\hfill  12\hfill 16 
& & 5\hfill  10\hfill  15\hfill  20\hfill 25\\ 
\hline
$z_i$ & $0$, $-4.26$, $-7.73$, 
$0. + 1.\times10^{50}i$, \quad
$0. - 1.\times10^{50}i$
& 
& $0$, $-5.52$, $-10.$, $-14.47$,\quad
$0. + 1.\times10^{50}i$,
$0. - 1.\times10^{50}i$
&
& $0$, $-6.77$, $-12.28$, $-17.71$, $-23.22$, \quad
$0. + 1.\times10^{50}i$,
$0. - 1.\times10^{50}i$
\\ 
\hline
\end{tabular}
\end{center}
\caption{Eigenvalues for $\lambda=10^{100}$}
\label{table:2}
\end{table}

\section{The Burgers model}\label{se:BM}
\begin{figure} [h!]
\centering
\begin{tikzpicture}[line cap=round,line join=round,>=triangle 45,x=0.5cm,y=0.5cm]
\clip(-0.13,-0.86) rectangle (12.62,5.29);
\fill[fill=black,pattern=north east lines] (0.25,3.5) -- (1,3.5) -- (1,0.5) -- (0.25,0.5) -- cycle;
\fill[fill=black,fill opacity=0.1] (7.25,3.5) -- (7.25,0.5) -- (7.37,0.5) -- (7.37,3.5) -- cycle;
\fill[fill=black,fill opacity=0.1] (10.87,3.51) -- (10.87,0.51) -- (10.99,0.51) -- (10.99,3.51) -- cycle;
\fill[fill=black,fill opacity=1.0] (12.25,2.08) -- (12.25,1.92) -- (12.4,2) -- cycle;
\draw [line width=1.2pt] (1,3.5)-- (1,0.5);
\draw [line width=1.2pt] (1,2)-- (1.74,1.99);
\draw [line width=1.2pt] (1.74,1.99)-- (1.89,2.48);
\draw [line width=1.2pt] (1.89,2.48)-- (2.19,1.51);
\draw [line width=1.2pt] (2.19,1.51)-- (2.48,2.48);
\draw [line width=1.2pt] (2.48,2.48)-- (2.78,1.51);
\draw [line width=1.2pt] (2.78,1.51)-- (3.06,2.48);
\draw [line width=1.2pt] (3.06,2.48)-- (3.36,1.52);
\draw [line width=1.2pt] (3.36,1.52)-- (3.5,2);
\draw [line width=1.2pt] (3.5,2)-- (5.25,2);
\draw [line width=3.6pt] (5.25,2.38)-- (5.25,1.62);
\draw [line width=1.2pt] (4.5,2.5)-- (6,2.5);
\draw [line width=1.2pt] (6,2.5)-- (6,1.5);
\draw [line width=1.2pt] (6,1.5)-- (4.5,1.5);
\draw [line width=1.2pt] (6,2)-- (7.25,2);
\draw [line width=0.4pt] (7.25,3.5)-- (7.25,0.5);
\draw (7.25,0.5)-- (7.37,0.5);
\draw (7.37,0.5)-- (7.37,3.5);
\draw (7.37,3.5)-- (7.25,3.5);
\draw [line width=1.2pt] (8.25,3.25)-- (8.4,3.75);
\draw [line width=1.2pt] (8.4,3.75)-- (8.69,2.75);
\draw [line width=1.2pt] (8.69,2.75)-- (8.98,3.75);
\draw [line width=1.2pt] (8.98,3.75)-- (9.27,2.75);
\draw [line width=1.2pt] (9.27,2.75)-- (9.56,3.75);
\draw [line width=1.2pt] (9.56,3.75)-- (9.86,2.75);
\draw [line width=1.2pt] (9.86,2.75)-- (10,3.25);
\draw [line width=3.6pt] (9.26,1.14)-- (9.26,0.36);
\draw [line width=1.2pt] (8.51,1.25)-- (10.01,1.25);
\draw [line width=1.2pt] (10.01,1.25)-- (10.01,0.25);
\draw [line width=1.2pt] (10.01,0.25)-- (8.51,0.25);
\draw [line width=1.2pt] (8.25,3.25)-- (7.37,3.25);
\draw [line width=1.2pt] (9.26,0.75)-- (7.37,0.76);
\draw (10.87,3.51)-- (10.87,0.51);
\draw (10.87,0.51)-- (10.99,0.51);
\draw (10.99,0.51)-- (10.99,3.51);
\draw (10.99,3.51)-- (10.87,3.51);
\draw [line width=1.2pt] (10,3.25)-- (10.87,3.25);
\draw [line width=1.2pt] (10.01,0.75)-- (10.87,0.75);
\draw [line width=1.2pt] (10.99,2)-- (12.25,2);
\draw [line width=1.2pt] (12.25,2.08)-- (12.25,1.92);
\draw [line width=1.2pt] (12.25,1.92)-- (12.4,2);
\draw [line width=1.2pt] (12.4,2)-- (12.25,2.08);
\draw (2.17,3.8) node[anchor=north west] {$E_1$};
\draw (4.9,3.8) node[anchor=north west] {$\eta_1$};
\draw (8.63,5) node[anchor=north west] {$E_2$};
\draw (8.63,0.16) node[anchor=north west] {$\eta_2$};
\end{tikzpicture}
\caption{Burgers model}\label{fig:Burgers}
\end{figure}
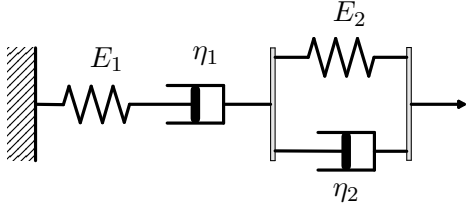

For reader's convenience,  first we will describe the Burgers Model, see e.g. \cite{SG}. The Maxwell and the Kelvin-Voigt two-element uniaxial models may be
investigated in this context and can be described by means of spring-dashpot systems. Indeed the Maxwell model consists of a linear elastic spring and a linear viscous dashpot element connected in a series, while the Kelvin-Voigt model is given by a linear spring element and a linear dashpot element which are connected in parallel. Those models are very simple, although they exhibit strong limitations.
In order to control such limitations, a more complex four-parameter (two YoungÕs modules $E_1$, $E_2$ and two viscosity parameters $\eta_1$, $\eta_2$) Burgers model which consists of two simple units, the Maxwell unit $(E_1,\eta_1 )$ and the Kelvin-Voigt unit $(E_2,\eta_2 )$ coupled in a series can be used, see Figure \ref{fig:Burgers}.

If the Burgers model is subject to a constant strain $\varepsilon=\varepsilon_0$ at $t=0$, a continuous stress relaxation modulus $E(t)$ is described by the combination of two exponential functions $e^{-r_1t}$ and $e^{-r_2t}$. Indeed, if we introduce the following parameters, whose definitions are due to Findley et al. \cite{FLO},
\begin{equation}\label{eq:Findley}
\begin{split}
&p_1=\frac{\eta_1}{E_1}+\frac{\eta_1}{E_2}+\frac{\eta_2}{E_2}
\qquad
p_2=\frac{\eta_1\eta_2}{E_1E_2}
\hskip1.5cm
q_1=\eta_1
\\
&q_2=\frac{\eta_1\eta_2}{E_2}
\hskip2.5cm
r_{1,2}=\frac{p_1\mp A}{2p_2}
\hskip1cm
A=\sqrt{p_1^2-4p_2}
\end{split}
\end{equation}
then the stress relaxation modulus is given by
\begin{equation}
E(t)=\frac{q_1-q_2r_1}{A}e^{-r_1t}-\frac{q_1-q_2r_2}{A}e^{-r_2t}
\,.
\end{equation}
From the definition of the stress relaxation modulus \eqref{const-rel} we can deduce that the integral kernel is given by $m(t)= -E'(t)$, so we have
\begin{equation}\label{eq:intm}
m(t)=\frac{r_1(q_1-q_2r_1)}{A}e^{-r_1t}-\frac{r_2(q_1-q_2r_2)}{A}e^{-r_2t}
\,.
\end{equation}
Because of \eqref{const-rel} and \eqref{eq:c2}, we have  $E(0)=\alpha+\int_0^\infty m(s)ds=c^2$. Since $E(0)=\int_0^\infty m(s)ds$ we have $\alpha=0$.
Set 
\begin{equation}\label{eq:def-b}
b_1:=\frac{r_1(q_1-q_2r_1)}{c^2A}\,,
\qquad
b_2:=-\frac{r_2(q_1-q_2r_2)}{c^2A}\,,
\end{equation}
and  note that
\begin{equation}
b_1>0,
\quad
b_2>0,
\quad
\frac{b_1}{r_1}+\frac{b_2}{r_2}=1
\,.
\end{equation}
Indeed, thanks to \eqref{eq:intm} and \eqref{eq:Findley} we get 
\begin{equation}
\frac{b_1}{r_1}+\frac{b_2}{r_2}=\frac1{c^2}\int_0^\infty m(s)ds=\frac{q_2}{c^2A}(r_2-r_1)=\frac{E(0)}{c^2}=1\,,
\end{equation}
\begin{equation}\label{eq:prop-c2}
c^2=\frac{q_2}{A}(r_2-r_1)=\frac{q_2}{p_2}
\,.
\end{equation}
Thanks to \eqref{eq:Findley} we note that
\begin{equation}
r_1+r_2=\frac{p_1}{p_2}
\,,
\qquad
r_1r_2=\frac{1}{p_2}
\,,
\end{equation}
and in view also of \eqref{eq:def-b} and \eqref{eq:prop-c2} we have
\begin{equation}\label{eq:}
\begin{split}
b_1+b_2
&=\frac{r_1q_1-q_2r_1^2-r_2q_1+q_2r_2^2}{c^2A}
=(r_1-r_2)\frac{q_1-q_2(r_1+r_2)}{c^2A}
\\
&=-\frac{1}{p_2^2}\frac{q_1p_2-q_2p_1}{c^2}
=\frac{p_1q_2-p_2q_1}{p_2q_2}
\,.
\end{split}
\end{equation}
Moreover
\begin{equation}\label{eq:accp}
r_1+ r_2-  b_1-  b_2=\frac{q_1}{q_2}
\,.
\end{equation}
In conclusion, the integro-differential equation \eqref{eq:problem0}  can be written in the form
\begin{equation}\label{eq:problem01}
u'' =c^2\triangle u-b_1c^2\int_0^t\ e^{-r_1(t-s)}\triangle u(s)ds-b_2c^2\int_0^t\ e^{-r_2(t-s)}\triangle u(s)ds
\,,
\end{equation}
where the constants $b_i$ and $r_i$ are defined in \eqref{eq:def-b} and \eqref{eq:Findley} respectively.
To give a complete spectral analysis of the above equation, we will transform it in a differential equation without integral terms.

%
First, we recast the integro-differential equation 
in an abstract setting. To this end, let 
$
H= L^2(\Omega)
$
be endowed with the usual scalar product and norm. 
We define the operator $L:D(L)\subset H\to H$  by
\begin{equation}\label{eq:operatorL}
\begin{array}{l}
D(L)=H^2(\Omega)\cap H_0^1(\Omega) \\
\\
L u=\displaystyle -c^2\triangle u\qquad u\in D(L)\,.
\end{array}
\end{equation}
It is well known that $L$ is a self-adjoint positive
 operator on $H$ with dense domain $D(L)$.
We denote by $\{\la_n\}_{n\ge1}$ a strictly increasing sequence  of  eigenvalues for the operator $L$ with
$\la_n>0$ and $\la_n\to\infty$ and we assume that the sequence of the corresponding eigenvectors $\{w_n\}_{n\ge1}$ constitutes a Hilbert basis for $H$.

Recalling that $b_i, r_i>0$, $i=1,2$, are defined in \eqref{eq:def-b} and \eqref{eq:Findley} and satisfy the condition $\frac{b_1}{r_1}+\frac{b_2}{r_2}=1$,
we consider the following Cauchy problem:
\begin{equation}\label{eq:system}
\begin{cases}
\displaystyle 
u''(t) +Lu(t)-b_1\int_0^t\ e^{-r_1(t-s)}L u(s)ds-b_2\int_0^t\ e^{-r_2(t-s)}L u(s)ds= 0
\hskip1cm  t\ge 0,
\\
u(0)=u_{0}\,,\quad u'(0)=u_{1}
\,.
\end{cases}
\end{equation}
For $(u_{0},u_{1})\in D(\sqrt{L})\times H$, we can write an expansion in terms of the eigenvectors $w_n$ of  the following type
\begin{equation}\label{eq:v0}
\begin{split}
& u_{0}=\sum_{n=1}^{\infty}u_{0n}w_{n}\,,\qquad\quad u_{0n}=
\langle u_{0},w_n\rangle \,,
\quad
\|u_0\|^2_{D(\sqrt{L})}=\sum_{n=1}^{\infty}u_{0n}^2\lambda_n\,,
\\
& u_{1}=\sum_{n=1}^{\infty}u_{1n}w_{n}\,,\qquad\quad u_{1n}=\langle u_{1},w_n\rangle \,,\quad
\|u_1\|^2_{H}=\sum_{n=1}^{\infty}u_{1n}^2\,.
\end{split} 
\end{equation} 
To write  the solution $u(t)$ of  \eqref{eq:system} as a series, that is
\begin{equation}\label{eq:rep0}
u(t)=\sum_{n=1}^{\infty}u_{n}(t)w_{n}\,,
\qquad
u_{n}(t)=\langle u(t),w_n\rangle
\,,
\end{equation}
we put that  expression for $u$ 
into \eqref{eq:system} and multiply  by $w_n$. 
It follows that  for any $n\in\N$ $u_{n}$ is the solution of the Cauchy problem
\begin{equation}\label{eq:secondsys}
\begin{cases}\displaystyle
u_{n}''
+\la_{n}u_{n}-\la_{n}b_1\int_0^t e^{-r_1(t-s)}u_{n}(s) ds-\la_{n}b_2\int_0^t e^{-r_2(t-s)}u_{n}(s) ds=0,
\\
u_{n}(0)=u_{0n}\,, \quad u_{n}'(0)=u_{1n}\,.
\end{cases}
\end{equation}
For a while, to simplify the notations we will drop the dependence on index $n$.
By means of derivations and integrations by parts one can establish that a scalar function $u$  defined on the interval $[0,\infty)$ is a solution of the second-order integro-differential equation 
\begin{equation}\label{eq:system0}
u{''}+\la u-\la b_1 \int_0^t e^{- r_1(t-s)}u(s) ds-\la b_2 \int_0^t e^{- r_2(t-s)}u(s) ds=0\,,
\qquad t\ge 0,
\end{equation}
if and only if $u$ is  a solution of the fourth-order differential equation 
\begin{equation}\label{eq:fourth}
\displaystyle 
u^{(4)}+( r_1+ r_2) u{'''}+(\la+ r_1 r_2) u{''}+\la ( r_1+ r_2-  b_1-  b_2) u'
=0,
\quad t\ge 0,
\end{equation}
and the conditions 
\begin{equation}\label{eq:fourth1}
\displaystyle
u{''}(0)=-\la u{}(0)
\,,
\qquad
u{'''}(0)=\la ( b_1+ b_2) u(0)-\la u{'}(0)
\end{equation}
are satisfied.
Therefore,
we have to evaluate the solutions of the $4^{\rm th}$--degree characteristic equation in the variable $z$
\begin{equation*}
z^{4}+(r_1+r_2) z^{3}+(\la_{n}+r_1r_2) z^{2}+\la_{n} (r_1+r_2- b_1- b_2)z
=0\,.
\end{equation*}
We have the solution $z=0$. To obtain the others, 
we have to solve the cubic equation
\begin{equation}\label{eq:fchar}
z^{3}+(r_1+r_2) z^{2}+(\la_{n}+r_1r_2) z+\la_{n} (r_1+r_2- b_1- b_2) 
=0\,.
\end{equation}
By means of the Cardano formula we have  the three solutions of \eqref{eq:fchar}: one is a real number $\rho_{n}$ and the others $i\omega_n$, $-i\overline{\omega_n}$ are  complex conjugate numbers. Moreover, $\rho_{n}$ and $\omega_n$ exhibit the following
asymptotic behavior as $n$ tends to $\infty$: 
\begin{multline}\label{eq:lambda1}
\rho_{n}=b_1+b_2-r_1-r_2
-{(b_1+b_2-r_1)(b_1+b_2-r_2)\over\la_{n}}(b_1+b_2-r_1-r_2)+O\Big({1\over{\la_{n}^{2}}}\Big)
\\
=b_1+b_2-r_1-r_2
+O\Big({1\over{\la_{n}}}\Big)
\,,
\end{multline}
\begin{multline}\label{eq:lambda2}
\omega_{n}=
\sqrt{\la_{n}}+{1\over8}\Big((b_1+b_2)(3(b_1+b_2)-4r_1)-4(b_1+b_2-r_1)r_2\Big){1\over\sqrt{\lambda_{n}}}
\\
+i \Big[{b_1+b_2\over 2}
-{(b_1+b_2-r_1)(b_1+b_2-r_2)\over2\la_{n}}(b_1+b_2-r_1-r_2) \Big]
+O\Big({1\over{\la_{n}^{3/2}}}\Big)
\\
=
\sqrt{\la_{n}}+i{b_1+b_2\over 2}
+O\Big({1\over{\sqrt{\la_{n}}}}\Big)
\,.
\end{multline}
We observe that for $b_2=r_2=0$ we have 
\begin{equation*}
\rho_{n}=b_1-r_1
-{b_1\big(b_1-r_1\big)^2\over\la_{n}}+O\Big({1\over{\la_{n}^{2}}}\Big),
\end{equation*}
\begin{equation*}
\omega_{n}=
\sqrt{\la_{n}}+{b_1\over2}\Big({3\over4}b_1-r_1\Big){1\over\sqrt{\lambda_{n}}}
+i \Big[{b_1\over 2}
-{b_1\big(b_1-r_1\big)^2\over2\la_{n}} \Big]
+O\Big({1\over{\la_{n}^{3/2}}}\Big)
\,,
\end{equation*}
that is the case of an integro-differential equation with a single exponential kernel given by $b_1e^{-r_1 t}$, see \cite{LoretiSforza, LoretiSforza1}.

Finally, taking also into account the conditions \eqref{eq:fourth1} and reintroducing the dependence on index $n$, $u_{n}$ is the solution of problem \eqref{eq:secondsys} if and only if  $u_{n}$ is the solution of the Cauchy problem
\begin{equation}\label{eq:third}
\begin{cases}
\displaystyle 
u_{n}^{(4)}+(r_1+r_2) u_{n}'''+(\la_{n}+r_1r_2) u_{n}''+\la_{n} (r_1+r_2- b_1- b_2) u_{n}'=0
\\
 u_{n}(0)=u_{0n},
\quad
 u_{n}'(0)=u_{1n},
 \\
 u_{n}''(0)=-\la_{n}u_{0n},
 \quad
 u_{n}'''(0)=\la_{n} (b_1+b_2) u_{0n}-\la_{n} u_{1n}
 \,.
\end{cases}
\end{equation} 
We are able to write the solution  $u_{n}(t)$ of (\ref{eq:third}) in the form
\begin{equation}\label{eq:f1j}
u_{n}(t)=R_{1,n}+R_{2,n}e^{\rho_{n} t}+
C_{n}e^{i\omega_{n} t}+\overline{C_{n}}e^{-i\overline{\omega_{n}}t}\,,
\end{equation}
where the coefficients $R_{1,n},R_{2,n}\in\R$ and $C_{n}\in\C$ can be determined by imposing
 the initial conditions. Therefore
we have to  solve 
 the system
\begin{equation}\label{vandermonde}
\left \{\begin{array}{l}
R_{1,n}+R_{2,n}+C_{n}+\overline{C_{n}}=u_{0n},\\ 
\\
\rho_{n} R_{2,n}+i\omega_{n} C_{n}-i\overline{\omega_{n}C_{n}}=u_{1n},\\
\\
\rho_{n}^2R_{2,n}-\omega_{n}^2C_{n}-\overline{\omega_{n}^2C_{n}}=-\la_{n}u_{0n},\\
\\
\rho_{n}^3R_{2,n}-i \omega_{n}^3C_{n}+i\overline{\omega_{n} ^3C_{n}}=\la_{n} (b_1+b_2) u_{0n}-\la_{n} u_{1n}.
\end{array}\right .
\end{equation}
Indeed, we obtain that the coefficients 
have the following asymptotic behavior as $n$ tends to $\infty$:
\begin{equation}\label{eq:asy_R1}
R_{1,n}={r_1r_2u_{1n}\over(r_1+r_2-b_1-b_2)\la_{n}}
+( u_{0n}+u_{1n})O\Big({1\over{\la_n^{2}}}\Big),
\end{equation}
\begin{equation}\label{eq:asy_R2}
R_{2,n}={(b_1+b_2-r_1)(b_1+b_2-r_2)\big(u_{0n}(b_1+b_2-r_1-r_2)+u_{1n}\big)\over(b_1+b_2-r_1-r_2)\la_{n}}
+( u_{0n}+u_{1n})O\Big({1\over{\la_n^{2}}}\Big),
\end{equation}
\begin{multline}\label{eq:asy_C}
C_{n}={u_{0n}\over 2}
-\frac{i}{4}
\big((b_1+b_2) u_{0n}+2u_{1n}\big)\frac{1}{\sqrt{\lambda_n}}
-\big((b_1+b_2-r_1)(b_1+b_2-r_2)u_{0n}
+(b_1+b_2)  u_{1n}\big)\frac1{2\la_{n}}
\\
+(u_{0n}+u_{1n})O\Big({1\over{\la_{n}^{3/2}}}\Big)
\,.
\end{multline}
Again we note that for $b_2=r_2=0$ we gain the result available for a single exponential kernel $b_1e^{-r_1 t}$, see \cite{LoretiSforza1}, that is
\begin{equation*}
R_{1,n}=0,
\quad
R_{2,n}={b_1\over\la_n}(u_{0n}(b_1-r_1)+u_{1n})
+( u_{0n}+u_{1n})O\Big({1\over{\la_n^{2}}}\Big),
\end{equation*}
\begin{equation*}
C_{n}={u_{0n}\over 2}
-\frac{i}{4}
\big(b_1 u_{0n}+2u_{1n}
\big)
\frac{1}{\sqrt{\lambda_n}}
-\frac{b_1}{2}
\big((b_1-r_1)u_{0n}
+u_{1n}\big)\frac1{\la_{n}}
+(u_{0n}+u_{1n})O\Big({1\over{\la_{n}^{3/2}}}\Big)
\,.
\end{equation*}
To get an explicit expression for the eigenvalues and eigenvectors of the operator $L$ defined by \eqref{eq:operatorL}
we will use polar coordinates.
First, we introduce the set ${\mathcal D}:=\{(r,\theta): 0<r<R\,, \ \theta\in[0,2\pi]\}$ and consider the operator $L$ in the space
$
H= L^2({\mathcal D})
$
endowed with the usual scalar product and norm 
\begin{equation*}
\langle u,v\rangle:=\int_0^R\int_0^{2\pi} ru(r,\theta)v(r,\theta)\ drd\theta\,,
\quad
\|u\|:=\left(\int_0^R\int_0^{2\pi} r|u(r,\theta)|^{2}\ drd\theta\right)^{1/2}
\quad
u,v\in L^2({\mathcal D})\,.
\end{equation*}
Moreover, we recall that  the Laplacian
in polar coordinates  is given by
\begin{equation*}
\triangle=\frac1{r}\frac{\partial}{\partial r}
\Big(r\frac{\partial}{\partial r}\Big)
+\frac1{r^2}\frac{\partial^2}{\partial\theta^2}
\,.
\end{equation*}
Therefore, we can rewrite the equation \eqref{eq:problem0} in the unknown $u(t,r,\theta)$
 \begin{equation}\label{eq:problem1}
u_{tt} =
\frac{c^2}{r}\big(ru_r\big)_r+\frac{c^2}{r^2}u_{\theta\theta}
-\frac1{r^2}\int_0^t\ m(t-s)
\Big(r\big(ru_r\big)_r+u_{\theta\theta}\Big)(s,r,\theta)ds
\quad t\ge 0\,,\,\, (r,\theta)\in {\mathcal D}.
\end{equation} 
For the sake of completeness, we briefly recall standard argumentations.
To determine the eigenvalues of the Laplacian, we have to solve
\begin{equation}\label{eq:eig}
-\triangle u(r,\theta)=\lambda^2 u(r,\theta)
\end{equation}
\begin{equation}
u(R,\theta)=0
\end{equation}
To this end, we attempt separation of variables by writing
\begin{equation*}
 u(r,\theta)=\Phi(r)\Theta(\theta).
\end{equation*}
Then \eqref{eq:eig} becomes
\begin{equation*}
r^2\frac{d^2 \Phi}{d r^2}\Theta+r\frac{d \Phi}{d r}\Theta
+\Phi\frac{d^2\Theta}{d\theta^2}  
+\lambda^2r^2 \Phi \Theta
=0
\,.
\end{equation*}
If we divide by $\Phi \Theta$, then we obtain
\begin{equation}\label{eq:eig1}
\frac{r^2}\Phi\frac{d^2 \Phi}{d r^2}
+\frac{r}{ \Phi}\frac{d \Phi}{d r}
+\frac1{\theta}\frac{d^2\Theta}{d\theta^2} 
+\lambda^2 r^2
=0
\,.
\end{equation}
The function $\Theta$ must be sinusoidal, that is
\begin{equation}\label{eq:theta}
\frac1{\theta}\frac{d^2\Theta}{d\theta^2} =-n^2,
\end{equation}
and hence, for $a_n\in\C$ we have
\begin{equation}\label{eq:thetas}
\Theta(\theta)=a_ne^{in\theta}+\overline{a_n}e^{-in\theta}.
\end{equation}
Plugging \eqref{eq:theta} into \eqref{eq:eig1}, we obtain
\begin{equation}\label{eq:R0}
r^2\frac{d^2 \Phi}{d r^2}
+r\frac{d \Phi}{d r}
+(\lambda^2 r^2-n^2)\Phi
=0
\,,
\end{equation}
with the boundary condition $\Phi(R)=0$.
We can eliminate $\lambda^2$ from the previous equation by making a change of variables. Indeed, if we set $x=\lambda r$, then the equation \eqref{eq:R0} becomes
\begin{equation}\label{eq:R1}
x^2\frac{d^2 \Phi}{d x^2}
+x\frac{d \Phi}{d x}
+(x^2-n^2)\Phi
=0
\,,
\end{equation}
which is called {\it Bessel's equation of order} $n$. A solution of \eqref{eq:R1} is given by
\begin{equation}\label{eq:func_Bess}
J_n(x)=\sum_{h=0}^\infty
\frac{(-1)^h}{h!(h+n)!}\Big(\frac{x}{2}\Big)^{n+2h}\,,
\end{equation}
which is called the {\it Bessel function of the first kind of order} $n$. It follows that a solution of \eqref{eq:R0} is given by $J_n(x)=J_n(\lambda r)$. The boundary condition $\Phi(R)=0$ is satisfied if
$
J_n(\lambda R)=0,
$
that is 
$
\lambda =\frac{\lambda_{nk}}{R}\,,
$
where $\lambda_{nk}$, $k\in\N$, are the positive zeros of $J_n$.
Therefore, the eigenvalues for  $L$ given by \eqref{eq:operatorL} are
 $c^2\big(\frac{\lambda_{nk}}{R}\big)^2$ and the corresponding eigenfunctions are $J_n\big(\frac{\lambda_{nk}}{R}r\big)e^{\pm in\theta}$, which form an orthogonal basis for $L^2({\mathcal D})$.

In order to simplify notations, we will define $J_{-n}$ to be the same as $J_n$ whenever $n$ is an integer:
\begin{equation*}
J_{-n}:=J_{n},
\quad
\lambda_{-nk}:=\lambda_{nk}
\,,\qquad
n\in\N\cup\{0\},\ \ k\in\N\,.
\end{equation*}
In conclusion, thanks to \eqref{eq:rep0}, \eqref{eq:f1j}, \eqref{eq:lambda1} and \eqref{eq:lambda2} we have the following representation for the solution
\begin{equation*}
u(t,r,\theta)
=
\sum_{n=-\infty}^{\infty}\sum_{k=1}^{\infty}
\Big(R_{1,nk}e^{in\theta}+R_{2,nk}e^{\rho_{nk}t+in\theta}+C_{nk}e^{i(\omega_{nk}t+n\theta)}
+\overline{C_{nk}}e^{-i(\overline{\omega_{nk}}t+n\theta)}\Big)
{J_{n}}\Big(\frac{\lambda_{nk}}{R} r\Big)\,,
\end{equation*}
where  
\begin{equation*}
\rho_{nk}=b_1+b_2-r_1-r_2
+O\Big({1\over{\la_{nk}^2}}\Big)
\,,
\quad
\omega_{nk}=
\frac{c}R\la_{nk}+i{b_1+b_2\over 2}
+O\Big({1\over{\la_{nk}}}\Big)
\,,
\end{equation*}
being $\lambda_{nk}$ the positive zeros of the Bessel function $J_n$ defined by \eqref{eq:func_Bess}.
\section{Conclusions}
In this paper we have investigated glass relaxation models, starting by a well-known model in literature, see e.g. \cite{MM} and references therein. Due to the complexity of the problem we have approximated the stretched exponential relaxation by means of a Prony series. For a general Prony series  we have established some partial results concerning the spectral analysis of the problem. In particular, by induction on the number of the terms of the Prony series the integro-differential equation showing the viscoelastic properties of the glass relaxation has always a null eigenvalue and the sum of all its eigenvalues is given by minus the sum of 
the exponents of the Prony series.

In order to give more accurate results, we simplified the problem by taking under consideration the Burgers model, where the Prony series consists of two decreasing exponential functions.
In that case we  have been able to
give a complete description of the oscillations of the material in its relaxation stage, when it shows viscoelastic features.
In particular, our analysis has revealed that the accumulation point of the branch of the real eigenvalues $\rho_n$, see \eqref{eq:lambda1},  depends only on the Kelvin-Voigt unit $(E_2,\eta_2 )$, see Figure \ref{fig:Burgers}. Indeed, taking into account \eqref{eq:accp} and \eqref{eq:Findley},
we have obtained
\begin{equation*}
b_1+b_2-r_1-r_2=-\frac{q_1}{q_2}=-\frac{E_2}{\eta_2}
\,.
\end{equation*}

\end{document}